\documentclass[a4paper]{amsart}
\usepackage{amsthm,amsfonts,amsmath,amssymb}
\usepackage[abs]{overpic}
\usepackage{thm-restate}

\newtheorem{theorem}{Theorem}[section]
\newtheorem{proposition}[theorem]{Proposition}

\newtheorem{remark}[theorem]{Remark}

\begin{document}
\title[The Gaps between the $\Delta$-Unknotting Number and its lower bounds]{The Gaps between the $\Delta$-Unknotting Number and its lower bounds}
\author{Kazumichi Nakamura}
\email{mi-ka-07130703@ezweb.ne.jp}





\begin{abstract}
Let $u^{\Delta}(K)$ denote the $\Delta$-unknotting number of a knot $K$, and let $a_2(K)$ denote the second coefficient of the Conway polynomial. It is known that
$u^{\Delta}(K)\geq |a_2(K)|$.
We define the $\Delta$-gap of $K$ by
$G^{\Delta}(K):=u^{\Delta}(K)-|a_2(K)|$.
In this paper, we investigate the behavior of this invariant for various families of knots.
For example, we prove that the $\Delta$-gap of the pretzel knot $K=P(2s, 4s-1,4s+1)$, where $s$ is a positive integer, satisfies $G^{\Delta}(K) \geq 4s$.
Consequently, the $\Delta$-gap is unbounded.
To further investigate the structure of the $\Delta$-gap, we distinguish between $\Delta$-moves that increase $a_2(K)$ and those that decrease it, since each $\Delta$-move changes $a_2(K)$ by $\pm1$. We analyze how these two types of moves contribute to the $\Delta$-gap.
We also investigate the gaps between $u^{\Delta}(K)$ and other knot invariants that provide lower bounds for it.
In particular, we establish explicit formulas and lower bounds for these gaps for certain families of pretzel knots, torus knots, two-bridge knots, and Montesinos knots. Our results show that these gaps can be arbitrarily large.
In addition, we identify prime positive knots with $u^{\Delta}(K)=2$ or $3$, and
propose a conjecture on a lower bound for $u^{\Delta}(K)$ for prime positive knots.
\end{abstract}

\maketitle

\noindent\textbf{keywords:}$\Delta$-move, $+\Delta$-move, $-\Delta$-move, $\Delta$-Gordian distance, $\Delta$-unknotting number, Conway polynomial, signature, 4-dimensional genus, genus, X-unknotting number, pretzel knot, torus knot, two-bridge knot, Montesinos knot, $\Delta$-gap.

\section{Introduction}
In this paper, we study the gaps between the 
$\Delta$-unknotting number and various lower bounds for it, including the absolute value of the second coefficient of the Conway polynomial.

\begin{figure}[htbp]
 \centering    \includegraphics[width=0.8\linewidth]
    {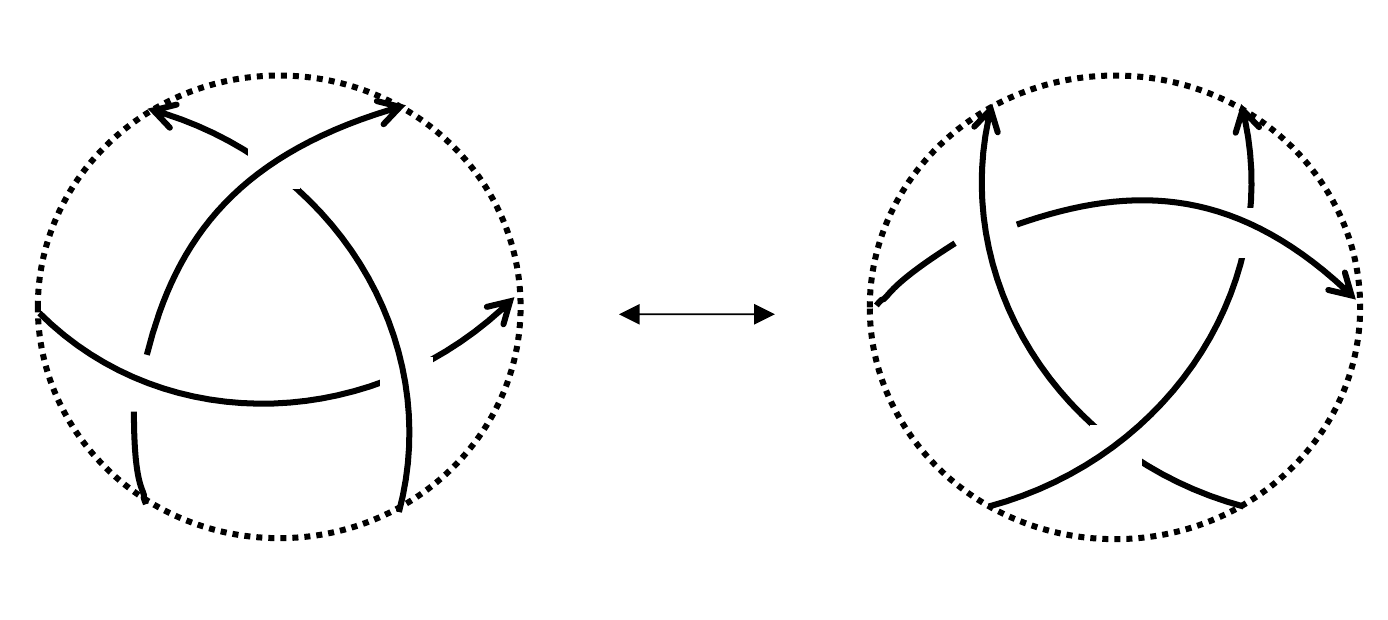}
 \caption{A $\Delta$-move.}
 \label{fig:delta}
\end{figure}

In \cite{MaN}, H. Murakami and Y. Nakanishi introduced a local move on regular diagrams of oriented knots and links, called a $\Delta$-move (or $\Delta$-unknotting operation), as illustrated in Figure \ref{fig:delta}.
They proved that two knots can be deformed into each other by a finite sequence of $\Delta$-moves. 
The $\Delta$-Gordian distance $d_G^{\Delta}(K, K^{'})$ of two oriented knots $K$ and $K^{'}$ is defined as the minimum number of $\Delta$-moves needed to deform a diagram of $K$ into that of $K{'}$. 
The $\Delta$-unknotting number $u^{\Delta}(K)$  of an oriented knot $K$ is defined as the $\Delta$-Gordian distance $d_G^{\Delta}(K, O)$ of $K$ and the trivial knot $O$.
In \cite{Oka1}, M. Okada proved that $u^{\Delta}(K) \geq |a_2(K)|$, where $a_2(K)$ denotes the second coefficient of the Conway polynomial of $K$. 

This inequality \(u^{\Delta}(K)\geq |a_2(K)|\) is quite effective: among the 249 knots with at most 10 crossings, equality holds for at least 196 knots.
This observation naturally leads us to consider the gap between $u^{\Delta}(K)$ and $|a_2(K)|$.

\begin{restatable}{definition}{defgap}\label{def:gap}
For a knot $K$, define
\[
G^{\Delta}(K) := u^{\Delta}(K) - |a_2(K)|.
\]
We call $G^{\Delta}(K)$ the $\Delta$-gap of $K$.
\end{restatable}
The $\Delta$-gap measures the gap between $u^{\Delta}(K)$ and the lower bound $|a_2(K)|$. It is natural to ask which values the $\Delta$-gap can take and how it behaves for various families of knots. Computational examples reveal several interesting behaviors.

\vspace{1em}

For example, $G^{\Delta}(K)=0$ for torus knots, a certain family of Turk's head knots, a certain family of two-bridge knots, positive pretzel knots, and a certain family of Montesinos knots.
On the other hand, there are many knots for which $G^{\Delta}(K)=2$. 
In \cite{Mura}, examples of knots satisfying $u^{\Delta}(K)=2$ and $a_2(K)=0$ are given.
These examples illustrate that the $\Delta$-gap contains information beyond the lower bound $u^{\Delta}(K) \geq |a_2(K)|$.

\vspace{1em}
In this paper, we investigate the $\Delta$-gap for a certain pretzel knots and a certain two-bridge knots. Our main result provides a lower bound for the $\Delta$-gap of these knots.

\begin{restatable}{theorem}{main}\label{thm:main}
Let $K=P(2s,4s-1,4s+1)$ be a pretzel knot, where $s$ is a positive integer.
Then, we have
\[
G^{\Delta}(K) \geq 4s.
\]
\end{restatable}

\begin{restatable}{theorem}{mains}\label{thm:mains}
Let $K=P(2s, 4s+1, 4s+1)$ be a pretzel knot, where $s$ is a positive integer.
Then, we have
\[
G^{\Delta}(K) \geq 3s.
\]
\end{restatable}

\begin{restatable}{theorem}{maint}\label{thm:maint}
Let $K=C(2s+1, -4s-2)$ be a two-bridge knot, where $s$ is a positive integer.
Then, we have
\[
G^{\Delta}(K) \geq 2s.
\]
Here, $C(2s+1, -4s-2) \cong P(2s,1,4s+1)
\cong C(2s,1,4s+1)$.
\end{restatable}

These results show, in particular, that the $\Delta$-gap is unbounded. 
Indeed, by taking $s$ arbitrarily large, the $\Delta$-gap of these knots can be made arbitrarily large. 

Thus, although the lower bound $|a_2(K)|$ provides useful information about the $\Delta$-unknotting number, the gap between $u^{\Delta}(K)$ and this lower bound can be arbitrarily large.

\vspace{1em}

In Section \ref{sec3}, in order to study its structure, we distinguish $\Delta$-moves according to whether they increase or decrease $a_2(K)$, since each $\Delta$-move changes $a_2(K)$ by $\pm1$. In Remark \ref{remark:gap34}, we investigate how these two types of moves contribute to the $\Delta$-gap.

\vspace{1em}


In Section \ref{sec4}, we investigate the gaps between $u^{\Delta}(K)$ and several lower bounds, including half the absolute value of the signature $\frac{|\sigma (K)|}{2}$, 4-dimensional genus $g^*(K)$, (Seifert) genus $g(K)$,
and $X$-unknotting number $u^X(K)$. 
We establish explicit formulas and lower bounds for these gaps for several families of knots, including pretzel knots, torus knots, two-bridge knots, and Montesinos knots.
In particular, our results demonstrate that the gaps between $u^{\Delta}(K)$ and these lower bounds can be arbitrarily large (see Theorems \ref{thm:mainu}, \ref{thm:mainx}, \ref{thm:mainv}, and \ref{thm:mainw}).

\vspace{1em}
In Section \ref{sec5}, we identify prime positive knots with $u^{\Delta}(K)=2$ or $3$ (see Theorem \ref{thm:thmpkdtwo}), 
and propose the following conjecture, which suggests a new lower bound for $u^{\Delta}(K)$ for prime positive knots.

\begin{restatable}{conjecture}{conjmain} \label{conj:conjmain}
Let $K$ be a prime positive knot with 
$c(K) = n$, where $c(K)$ denotes the crossing number of $K$ and $n$ is a positive integer.
Then, we have
\[
u^{\Delta}(K)\geq
\begin{cases}
\frac{1}{2}(n-1), & \text{if $n$ is odd},\\[6pt]
\frac{1}{2}n, & \text{if $n$ is even}.
\end{cases}
\]

Moreover, these bounds are sharp. 

Indeed, for every positive odd
integer $n$, the knot
$K=C(n-2,2)=P(1,1,n-2)$
satisfies $c(K)=n$ and $u^{\Delta}(K)=\frac{1}{2}(n-1)$.

For every even integer $n\geq 8$, the knot
$K=M(0;(2,1),(3,2),(2n-13,2))$
satisfies $c(K)=n$ and $u^{\Delta}(K)=\frac{1}{2}n$.
\end{restatable}

Conjecture \ref{conj:conjmain} holds for prime positive knots with at most 10 crossings.
If it holds in general, this conjecture could serve as one of the criteria for determining whether a knot is positive.

\vspace{1em}
\section{Preliminaries}\label{sec2}


\subsection{Pretzel Knots} 
\label{subsec21}

\leavevmode

\begin{figure}[htbp]
    \centering
    \includegraphics[width=0.8\linewidth]{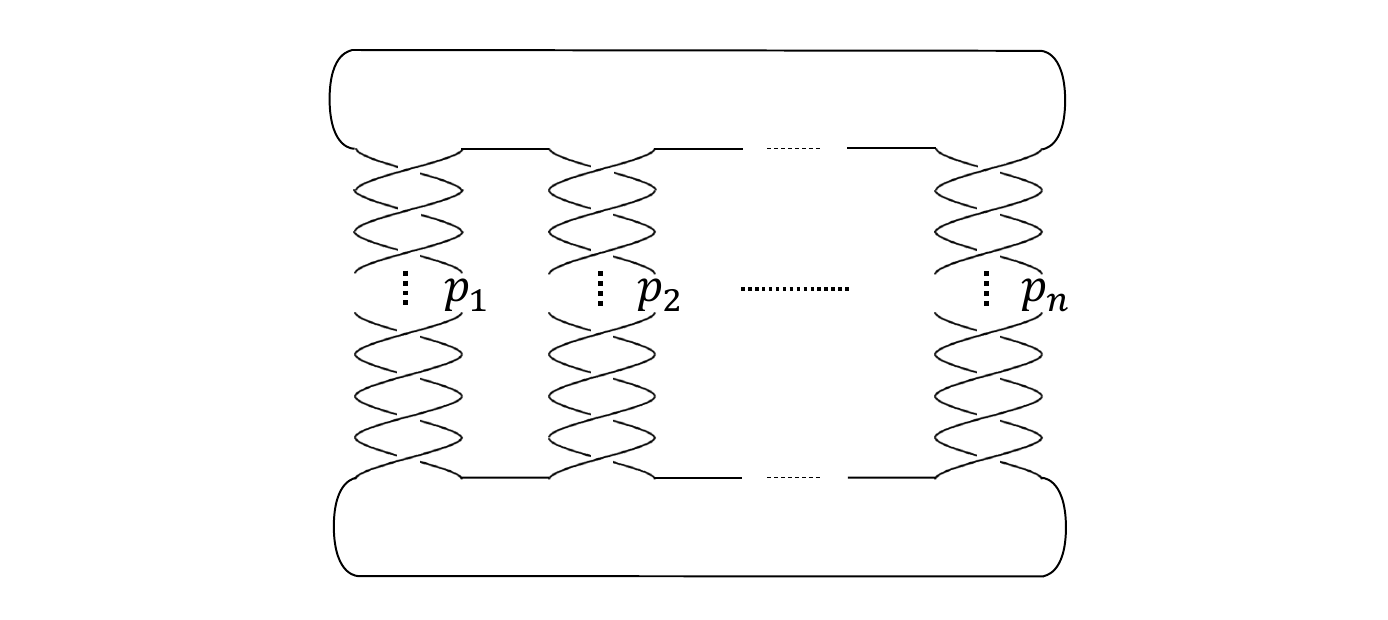}
    \caption{A pretzel knot $P(p_1, p_2, ... , p_n)$.}
    \label{fig:pretzel}
\end{figure}

A pretzel link $P(p_1, p_2, ... , p_n)$ is a link obtained by connecting $n$ twisted bands, where each band has $p_i$ half-twists.
The sign of $p_i$ determines the direction of the $i$th twist region
(see Figure \ref{fig:pretzel}).
A pretzel link is a knot if and only if precisely one $p_i$ is even, or all $p_i$ are odd and $n$ is odd.
Previous work \cite{Naka4} computed the $\Delta$-unknotting numbers of positive pretzel knots.

\subsection{Two-Bridge Knots} 
\label{subsec22}

\leavevmode

\begin{figure}[htbp]
 \centering
 \includegraphics[width=1\linewidth]{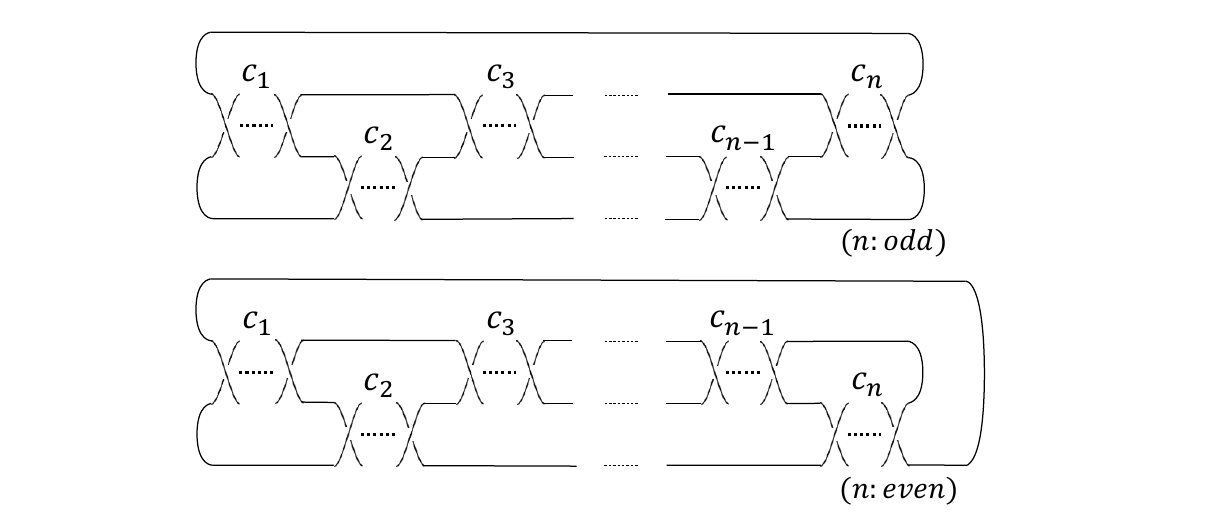}
 \caption{A two-bridge knot $C(c_1, c_2, ... , c_n)$.}
 \label{fig:twobridge}
\end{figure}

A knot $K$ is said to be a two-bridge knot if $K$ has a diagram as in Figure \ref{fig:twobridge}, called Conway's normal form. For a knot diagram as in Figure \ref{fig:twobridge}, each $|c_i|$ presents the number of half-twists for integers $c_1, c_2, ... , c_n$. In this paper, for the sign of $\alpha_i$, we assume that a right-handed half-twist is positive when $i$ is odd, and a left-handed half-twist is positive when $i$ is even. We denote this knot diagram by $C(c_1, c_2, ... , c_n)$.
Previous work \cite{Naka3} computed the $\Delta$-unknotting numbers of certain families of two-bridge knots.

\subsection{Montesinos Knots} 
\label{subsec23}

\leavevmode

\begin{figure}[htbp]
    \centering
    \includegraphics[width=0.8\linewidth]{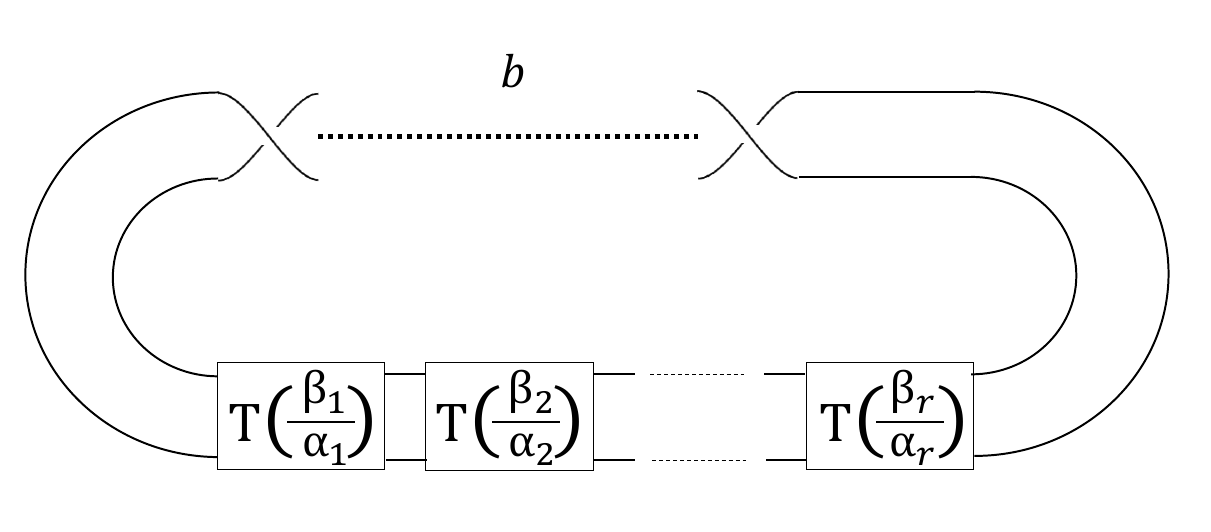}
    \caption{A Montesinos knot $M\big(b ; (\alpha_1, \beta_1), (\alpha_2, \beta_2), \ldots, (\alpha_r, \beta_r)\big)$.}
    \label{fig:montesinos}
\end{figure}

\begin{figure}[htbp]
    \centering
    \includegraphics[width=0.8\linewidth]{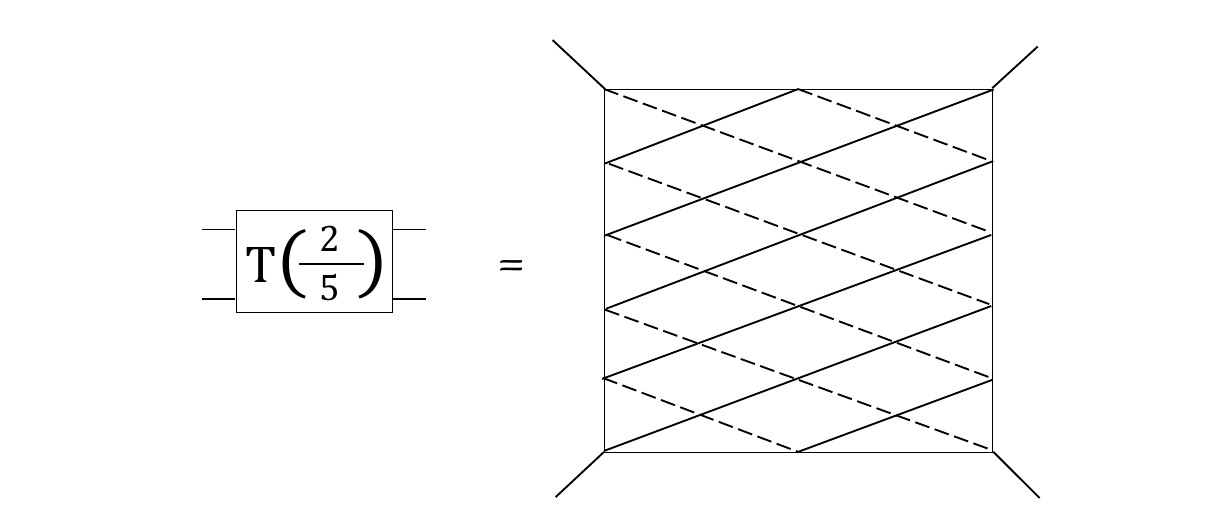}
    \caption{The tangle $T(2/5)$.}
    \label{fig:tangle}
\end{figure}

A Montesinos knot $K = M\big(b ; (\alpha_1, \beta_1), (\alpha_2, \beta_2), \ldots, (\alpha_r, \beta_r)\big)$ with $r$ branches is a knot as illustrated in Figure \ref{fig:montesinos}.
Here $r$, $b$, $\alpha_k$, and $\beta_k$ are integers 
such that $r \geq 0$,  $\alpha_k \geq 2$, and $\gcd(\alpha_k, \beta_k)=1$.
The quantity $|b|$ denotes the number of half-twists corresponding to the integer parameter $b$, where a right-handed half-twist is considered positive.
Let $T(\beta_k / \alpha_k)$ denote the rational tangle of slope $\beta_k / \alpha_k$ (see Figure  \ref{fig:tangle}).

\vspace{1em}
The upper left and lower left corners of $T(\beta / \alpha)$ are connected by a strand if and only if $\alpha$ is even.
Hence we see that if $K = M\big(0 ; (\alpha_1, \beta_1), (\alpha_2, \beta_2), \ldots, (\alpha_r, \beta_r)\big)$
is to be a knot rather than a link, at most one of $\alpha_1, \alpha_2$, ... , $\alpha_r$ can be even.
Since a cyclic permutation indices does not change the knot type, we hereafter assume that $\alpha_2$, ..., $\alpha_r$ are odd.
With this convention, we say that $K$ is of odd type if $\alpha_1$ is odd, and of even type if $\alpha_1$ is even.
Previous work \cite{Naka5} computed the $\Delta$-unknotting numbers of certain families of Montesinos knots.

\vspace{1em}
Since $a_2(K)=a_2(K^*)$ and $u^{\Delta}(K)=u^{\Delta}(K^*)$, where $K^*$ denotes the mirror image of a knot $K$, we do not distinguish a knot from its mirror image in this paper.

\vspace{1em}
\section{The ${\Delta}$-gap of knots}\label{sec3}

\subsection{Definition of the $\Delta$-gap} 
\label{subsec31}

\leavevmode
\vspace{1em}

The following proposition gives a useful relation between the $\Delta$-unknotting number and the second coefficient of the Conway polynomial.

\begin{proposition}[\cite{Oka1}]\label{prop:a2}
For any two knots $K$ and $K^{'}$, the difference $d_G^{\Delta}(K,K^{'}) - |a_2(K) - a_2(K^{'})|$ is a non-negative even integer. 
In particular, the difference $u^{\Delta}(K) -|a_2(K)|$ is also a non-negative even integer.
\end{proposition}

This proposition naturally leads us to consider the gap between $u^{\Delta}(K)$ and $|a_2(K)|$.

\vspace{1em}
\defgap*

Thus, $G^{\Delta}(K)$ is a knot invariant taking values in the non-negative even integers.

\vspace{1em}
\subsection{Examples of the $\Delta$-gap} 
\label{subsec32}

\leavevmode

\vspace{1em}

For example, $G^{\Delta}(O)
= u^{\Delta}(O) - |a_2(O)| = 0$.

Among the 249 knots with at most 10 crossings, 
$u^{\Delta}(K)=|a_2(K)|$ (i.e., $G^{\Delta}(K)=0$) holds for at least 196 knots.


\begin{table}[ht]
\centering
\small
\renewcommand{\arraystretch}{1.25}
\caption{$G^{\Delta}(K)$ among knots with at most 10 crossings.}
\label{tab:GDelta}
\begin{tabular}{c|rrrrrrrr|r}
\hline
Crossing number
& 3 & 4 & 5 & 6 & 7 & 8 & 9 & 10 & Total \\
\hline
Number of knots
& 1 & 1 & 2 & 3 & 7 & 21 & 49 & 165 & 249 \\
\hline
$G^{\Delta}(K)=0$
& 1 & 1 & 2 & 3 & 7 & 17 & 41 & 124 & 196 \\
$G^{\Delta}(K)=2$
& 0 & 0 & 0 & 0 & 0 & 4 & 7 & 25 & 36 \\
$G^{\Delta}(K)=4$
& 0 & 0 & 0 & 0 & 0 & 0 & 0 & 2 & 2 \\
\hline
$G^{\Delta}(K)=0$ or $2$
(undetermined)
& 0 & 0 & 0 & 0 & 0 & 0 & 1 & 12 & 13 \\
$G^{\Delta}(K)=2$ or $4$
(undetermined)
& 0 & 0 & 0 & 0 & 0 & 0 & 0 & 2 & 2 \\
\hline
\end{tabular}
\end{table}

Furthermore, it is known that
$$
u^{\Delta}(K)=|a_2(K)|
\qquad
\text{(i.e., \(G^{\Delta}(K)=0\))}
$$

for the following families of knots:
\begin{itemize}
\item torus knots \cite{NaNaU};
\item a certain family of Turk's head knots \cite{NaYa};
\item a certain family of two-bridge knots \cite{Naka3};
\item positive pretzel knots \cite{Naka4};
\item a certain family of Montesinos knots \cite{Naka5}.
\end{itemize}


\vspace{1em}

The following proposition provides examples of knots with $G^{\Delta}(K)=0$ and $G^{\Delta}(K)=2$.

\begin{proposition}[\cite{Mura}]\label{prop:Mura}
Let
$f(z) = 1 + b_1z^2 + b_2z^4 + \cdots + b_nz^{2n}$
be a polynomial with \(b_i \in \mathbb{Z}\) \((i=1,2,\ldots,n)\). If \(b_1 \neq 0\), then there exists a knot \(K\) with
$u^{\Delta}(K) = |b_1|$ 
and
$\nabla_K(z) = f(z)$.
If \(b_1 = 0\), there exists a knot \(K\) with
$u^{\Delta}(K) = 2$
and
$\nabla_K(z) = f(z)$.
\end{proposition}

\vspace{1em}


Among prime knots with at most 10 crossings satisfying $G^{\Delta}(K)=4$, apart from $10_{56}$ and $10_{98}$, the only ones are $10_{46}$ and $10_{162}$.


\begin{table}[htbp]
\centering
\caption{The values of $u^{\Delta}(K)$ and $a_2(K)$.}
\label{tab:u-a2}
\begin{tabular}{c|cc}
\hline
$K$ & $u^{\Delta}(K)$ & $a_2(K)$ \\
\hline
$10_{46}$  & 4 & 0 \\
$10_{56}$  & 2 or 4 & 0 \\
$10_{98}$  & 2 or 4 & 0 \\
$10_{162}$ & 7 & 3 \\
\hline
\end{tabular}
\end{table}

\subsection{Main Theorems on the $\Delta$-gap} 
\label{subsec33}

\leavevmode

\vspace{1em}

As shown in the previous subsection, $G^{\Delta}(10_{46}) = 4$, which is relatively large compared with those of the other knots with at most 10 crossings.

Since $10_{46}$ is a three-strand pretzel knot, we extend this observation to the family of three-strand pretzel knots and obtain the following theorems.

In the proofs of the main theorems, we use the inequality $u^{\Delta}(K) \geq \frac{|\sigma(K)|}{2}$.
For the computation of $a_2(K)$, we use the result of \cite{Naka4}.

\main*

\begin{proof}
Let $K=P(2s, 4s-1, 4s+1)$, where $s$ is a positive integer.

Recall that for a pretzel knot $P(p, q, r)$, where $p$ is a positive even integer and $q,r$ are positive odd integers,
we have 
$\sigma(P(p, q, r)) = 2-q-r$.

Therefore,
\[
\sigma(K) = 2-8s.
\]

Moreover,
\[
a_2(K)=0.
\]

Thus,
\[
G^{\Delta}(K) =  u^{\Delta}(K) - |a_2(K)|
\geq \frac{|\sigma(K)|}{2} = \frac{|2-8s|}{2} =4s-1.
\]

Since $G^{\Delta}(K)$ is a positive even integer,
it follows that
\[
G^{\Delta}(K) \geq 4s.
\]

The proof is complete.
\end{proof}

Let $K=P(2,3,5)=10_{46}$.
Then, we have $u^{\Delta}(K)=4$ and $a_2(K)=0$.
Thus $G^{\Delta}(K)=4$.
Therefore, equality holds in the inequality $G^\Delta(K)\geq 4s$ in Theorem \ref{thm:main}.

\mains*

\begin{proof}
Let $K=P(2s, 4s+1, 4s+1)$, where $s$ is a positive integer.

Therefore,
\[
\sigma(K) = -8s.
\]

Moreover,
\[
a_2(K)=s.
\]

Thus,
\[
G^{\Delta}(K) =  u^{\Delta}(K) - |a_2(K)|
\geq \frac{|\sigma(K)|}{2} - |a_2(K)| = \frac{|-8s|}{2} - s = 3s.
\]

Since $G^{\Delta}(K)$ is a positive even integer,
it follows that
\[
G^{\Delta}(K) \geq 3s.
\]

The proof is complete.
\end{proof}


\vspace{1em}
\maint*

\begin{proof}
Let $K=P(2s, 1, 4s+1) = C(2s+1, -4s-2)$, where $s$ is a positive integer.

Therefore,
\[
\sigma(K) = -4s.
\]

Moreover,
\[
a_2(K)=0.
\]

Thus,
\[
G^{\Delta}(K) =  u^{\Delta}(K) - |a_2(K)|
\geq \frac{|\sigma(K)|}{2} = \frac{|-4s|}{2} = 2s.
\]

The proof is complete.
\end{proof}

Let $K=C(3,-6)=P(2,1,5)=8_{2}$.
Then, we have $u^{\Delta}(K)=2$ and $a_2(K)=0$.
Thus $G^{\Delta}(K)=2$.
Therefore, equality holds in the inequality $G^\Delta(K)\geq 2s$ in Theorem \ref{thm:maint}.

\vspace{1em}
\subsection{Remarks on the $\Delta$-gap} 
\label{subsec34}

\leavevmode

\vspace{1em}

Let $K'$ be obtained from $K$ by a single ${\Delta}$-move. Then we have $|a_2(K')-a_2(K)|=1$.
If $a_2(K')-a_2(K)=1$, $K'$ is said to be obtained from $K$ by a single $+{\Delta}$-move.
If $a_2(K')-a_2(K)=-1$, $K'$ is said to be obtained from $K$ by a single $-{\Delta}$-move.

In other words, a $+\Delta$-move is a $\Delta$-move that increases $a_2$, whereas a $-\Delta$-move is a $\Delta$-move that decreases $a_2$.

\begin{remark}[\cite{NaYa}]\label{remark:gap33}
If a knot $K$ is deformed into a trivial knot $O$ by a finite sequence of $+\Delta$-moves (or $-\Delta$-moves), then we have 
$u^{\Delta}(K) = |a_2(K)|$.
\end{remark}

$$
\begin{aligned}
K \quad (a_2(K)<0)
&\xrightarrow{\ +\Delta\ }
K_1
\xrightarrow{\ +\Delta\ }
K_2
\xrightarrow{\ +\Delta\ }
\cdots
\xrightarrow{\ +\Delta\ }
K_n=O,
\\[6pt]
\text{or} \qquad
K \quad (a_2(K)>0)
&\xrightarrow{\ -\Delta\ }
K_1
\xrightarrow{\ -\Delta\ }
K_2
\xrightarrow{\ -\Delta\ }
\cdots
\xrightarrow{\ -\Delta\ }
K_n=O.
\end{aligned}
$$

$$
\Longrightarrow\qquad
u^\Delta(K)=|a_2(K)|=n.
$$

\vspace{1em}

The following observation relates the $\Delta$-gap to the numbers of $+\Delta$-moves and $-\Delta$-moves in a minimal $\Delta$-unknotting sequence.

\vspace{1em}

$$
K
\xrightarrow{\ \pm\Delta\ }
K_1
\xrightarrow{\ \pm\Delta\ }
K_2
\xrightarrow{\ \pm\Delta\ }
\cdots
\xrightarrow{\ \pm\Delta\ }
K_{n-1}
\xrightarrow{\ \pm\Delta\ }
K_n=O,
$$

$$
\text{where each move is either a \(+\Delta\)-move or a \(-\Delta\)-move.}
$$

\begin{remark}\label{remark:gap34}
Let $t_{+\Delta}(K)$ and $t_{-\Delta}(K)$ denote the numbers of $+\Delta$-moves and $-\Delta$-moves, respectively, required to realize $u^{\Delta}(K)$.
Then,  
\[
u^{\Delta}(K) = t_{+\Delta}(K) + t_{-\Delta}(K),
\]
and
\[
a_2(K) = - t_{+\Delta}(K) + t_{-\Delta}(K).
\]
Moreover, we have the following cases:
\begin{enumerate}
\item
If $a_2(K) > 0$, then
\[
G^{\Delta}(K) = 2t_{+\Delta}(K) < 2t_{-\Delta}(K).
\]
\item
If $a_2(K) < 0$, then
\[
G^{\Delta}(K) = 2t_{-\Delta}(K) < 2t_{+\Delta}(K).
\]
\item
If $a_2(K) = 0$, then
\[
G^{\Delta}(K) = u^{\Delta}(K) = 2t_{+\Delta}(K) = 2t_{-\Delta}(K).
\]
\end{enumerate}
In particular, we obtain
\[
G^{\Delta}(K) = \min \{ 2t_{+\Delta}(K), \quad 2t_{-\Delta}(K) \}.
\]
\end{remark}

\vspace{1em}

Thus, $G^{\Delta}(K)$ is determined by the smaller of the numbers of $+\Delta$-moves and $-\Delta$-moves in a minimal $\Delta$-unknotting sequence.

\vspace{1em}
\section{The gaps between the $\Delta$-unknotting number and other knot invariants}\label{sec4}

It is known that $u^{\Delta}(K)$ can be bounded below by several
knot invariants (see \cite{MaN}).
In particular, we have
\[
\begin{gathered}
u^\Delta(K) \geq g^*(K) \geq \frac{|\sigma(K)|}{2},\\
u^\Delta(K) \geq \frac{u^X(K)}{2}.
\end{gathered}
\]
Here, $g^*(K)$ denotes 4-dimensional genus of a knot $K$,
$\sigma(K)$ denotes its signature, and $u^X(K)$ denotes
its $X$-unknotting number.

According to \cite{Sto}, for a positive knot, we have
\[
\begin{gathered}
u^\Delta(K) \geq a_2(K) \geq u^X(K) \geq g(K).
\end{gathered}
\]
Here, $g(K)$ denotes genus of a knot $K$.
It is known that $g(K)=g^*(K)$ 
for every positive knot $K$ (see \cite{NakaT}).

We obtain the following theorems on the gaps between $u^{\Delta}(K)$ and these lower bounds.
In the proofs of these theorems, we use the inequality $u^{\Delta}(K) \geq |a_2(K)|$. 
For the computation of $u^{\Delta}(K)$ and $a_2(K)$, we use the results of \cite{Naka4,Naka3}.

\begin{restatable}{theorem}{mainu}\label{thm:mainu}
Let $K=P(p, q, r)$ be a pretzel knot of odd type, where $p$,$q$, and $r$ are odd integers satisfying $pq + pr + qr < 0$.
Then, we have
\[
u^{\Delta}(K) - \frac{|\sigma (K)|}{2} \geq \frac{|pq + pr + qr + 1|}{4}.
\]

In particular, if $K=P(p, -p, r)$,
then
\[
u^{\Delta}(K) - \frac{|\sigma (K)|}{2} \geq \frac{p^2 -1}{4}.
\]
Here, $u^{\Delta}(K) = \frac{p^2 -1}{4}$ and $\frac{|\sigma (K)|}{2} = 0$.
\end{restatable}

\begin{proof}
Let $K = P(p,q,r)$ be a pretzel knot of odd type, where $p$,$q$, and $r$ are odd integers satisfying $pq + pr + qr <0$. 
By \cite{GaJ}, we have
\[
\sigma (P(p,q,r)) = 0.
\]

Moreover,
\[
a_2(P(p,q,r))=\frac{pq+pr+qr+1}{4}.
\]

Thus,
\[
u^{\Delta}(K) - \frac{|\sigma (K)|}{2} 
\geq |a_2(K)| = \frac{|pq+pr+qr+1|}{4}.
\]

In particular, for $P(p,-p,r)$, we obtain
\[
u^{\Delta}(P(p,-p,r)) - \frac{|\sigma (P(p,-p,r))|}{2} 
\geq \frac{|-p^2+1|}{4} = \frac{p^2-1}{4}.
\]

The proof is complete.
\end{proof}

\begin{restatable}{theorem}{mainx}\label{thm:mainx}
Let $K=T(p,q)$ be a torus knot, where $p$ and $q$ are positive integers with $gcd(p,q) = 1$.
Then, we have
\begin{align*}
u^{\Delta}(K) - \frac{|\sigma (K)|}{2} 
&= u^{\Delta}(K) - g^*(K) \\
&= u^{\Delta}(K) - g(K)   \\
&= u^{\Delta}(K) - u^X(K) \\
&= \frac{1}{24}(p-1)(q-1)(pq+p+q-11),
\end{align*}
and
\[
u^{\Delta}(K) - \frac{u^X(K)}{2} 
= \frac{1}{24}(p-1)(q-1)(pq+p+q-5).
\]
Here, $u^{\Delta}(K) = \frac{1}{24}(p^2 -1)(q^2 -1)$ and $\frac{|\sigma (K)|}{2} = g^*(K) = g(K) = u^X(K) = \frac{1}{2}(p -1)(q -1)$.

In particular, if $K=T(p, p+1)$,
then
\[
u^{\Delta}(K) - \frac{u^X(K)}{2} 
= \frac{1}{24}p(p-1)(p-2)(p+5).
\]
\end{restatable}

\begin{proof}
Let $K=T(p,q)$ be a torus knot, where $p$ and $q$ are positive integers with $gcd(p,q) = 1$.
By \cite{NaNaU}, we have $u^{\Delta}(K) = \frac{1}{24}(p^2 -1)(q^2 -1)$.
It is well known that $\frac{|\sigma (K)|}{2} = g^*(K) = g(K) = u^X(K) = \frac{1}{2}(p -1)(q -1)$.
Therefore, the desired result follows, proving Theorem \ref{thm:mainx}.
\end{proof}

\vspace{1em}

\begin{restatable}{theorem}{mainv}\label{thm:mainv}
Let $K=C(2s+1,2)$ be a two-bridge knot, where $s$ is a positive integer.
Then, we have
\begin{align*}
u^{\Delta}(K) - \frac{|\sigma (K)|}{2} 
&= u^{\Delta}(K) - g^*(K) \\
&= u^{\Delta}(K) - g(K)   \\
&= u^{\Delta}(K) - u^X(K) \\
&= s,
\end{align*}
and
\[
u^{\Delta}(K) - \frac{u^X(K)}{2} 
= s + \frac{1}{2}.
\]
Here, $u^{\Delta}(K) = s + 1$ and $\frac{|\sigma (K)|}{2}= g^*(K) = g(K) = u^X(K) = 1$.
\end{restatable}

\begin{proof}
Let $K = C(2s+1,2)$ be a two-bridge knot, where $s$ is a positive integer. 

Since $K$ can be deformed into a trivial knot by a single X-move, we have,
\[
u^X(K) = 1.
\]

Since $K$ is a positive knot, we have 
$g^*(K)=g(K)=1$ and $\frac{|\sigma (K)|}{2}=1$.

Moreover,
\[
u^{\Delta}(K)=|a_2(K)|=\frac{2(2+4s+2)}{8}=s+1.
\]

Consequently,
\[
u^{\Delta}(K) - \frac{|\sigma (K)|}{2} 
= s + 1 -1 = s,
\]

and

\[
u^{\Delta}(K) - \frac{u^X(K)}{2} 
= s + 1 - \frac{1}{2} = s + \frac{1}{2}.
\]

The proof is complete.
\end{proof}

\begin{restatable}{theorem}{mainw}\label{thm:mainw}
Let $K=M(0;(2,1),(3,2),(4s-5,2))$ be a Montesinos knot, where $s$ is a positive integer with $s \geq 2$.
Then, we have
\begin{align*}
u^{\Delta}(K) - \frac{|\sigma (K)|}{2} 
&= u^{\Delta}(K) - g^*(K) \\
&= u^{\Delta}(K) - g(K)   \\
&= u^{\Delta}(K) - u^X(K) \\
&= s,
\end{align*}
and
\[
u^{\Delta}(K) - \frac{u^X(K)}{2} 
= s + 1.
\]
Here, $u^{\Delta}(K) = s + 2$ and $\frac{|\sigma (K)|}{2} = g^*(K) = g(K) = u^X(K) = 2$.
\end{restatable}



\begin{proof}
Let $K = M(0;(2,1),(3,2),(4s-5,2))$ be a Montesinos knot, where $s$ is a positive integer with $s \geq 2$. 

Since $K$ can be deformed into a trivial knot by two X-moves, we have
\[
u^X(K) \leq 2.
\]

On the other hand, since $K$ is a positive knot, it follows from \cite{PaT} that
\[
\frac{|\sigma (K)|}{2} \geq 2.
\]

Combining this with the inequality
\[
u^X(K) \geq \frac{|\sigma (K)|}{2},
\] 
we obtain
\[
u^X(K) = \frac{|\sigma (K)|}{2} =2.
\] 

Furthermore, since $K$ is a positive knot, we have 
$g(K) \leq 2$ and $g^*(K) \geq 2$.
Consequently, $g^*(K)=g(K)=2$.

Moreover, by the same argument as
in Remark \ref{remark:remarkmonte}, we obtain
\[
u^{\Delta}(K)=|a_2(K)|=s+2.
\]

Therefore,
\[
u^{\Delta}(K) - \frac{|\sigma (K)|}{2} 
= s + 2 -2 = s,
\]

and

\[
u^{\Delta}(K) - \frac{u^X(K)}{2} 
= s + 2 - \frac{2}{2} = s + 1.
\]

The proof is complete.
\end{proof}

The above results demonstrate that the gaps between 
$u^{\Delta}(K)$ and the classical lower bounds can be arbitrarily large for certain families of pretzel knots, torus knots, two-bridge knots, and Montesinos knots.

\vspace{1em}
\section{A conjecture on the $\Delta$-unknotting number for prime positive knots}\label{sec5}


We recall the following results from \cite{Naka4} 
and \cite{Sto}.

\begin{restatable}[\cite{Naka4}]{proposition}{pkdone}
\label{thm:thmpkdone}
A positive knot has $\Delta$-unknotting number one if and only if it is a right-handed trefoil.
\end{restatable}

\begin{restatable}[\cite{Naka4}]{proposition}{pomin} \label{pro:pomin}
Among positive pretzel knots of odd type with a fixed odd crossing number $n$, the ${\Delta}$-unknotting number is minimized by $P(1, 1, n-2)$, and the minimum value is $\frac{1}{2}(n-1)$.
\end{restatable}

\begin{proposition}[\cite{Sto}]\label{prop:Sto}
In a reduced positive $c$ crossing diagram, 
$a_2(K) \geq c/4$.
\end{proposition}

Since $u^{\Delta}(K) \geq |a_2(K)|$ and $c \geq c(K)$, it follows that $u^{\Delta}(K) \geq c(K)/4$.
Therefore, for any fixed value of $u^{\Delta}(K)$,
there are only finitely many knots $K$ satisfying that value.

In connection with the above results, we now state the following Theorems \ref{thm:thmpkdtwo} and \ref{thm:deltatbmin}.

\begin{restatable}{theorem}{pkdtwo} \label{thm:thmpkdtwo}
The following statements hold:
\begin{enumerate}
\item 
A prime positive knot has $\Delta$-unknotting number two if and only if it is a right-handed $5_2$ knot.
\item 
A prime positive knot has $\Delta$-unknotting number three if and only if it is a right-handed $5_1$ or $7_2$ knot.
\end{enumerate}
\end{restatable}

\begin{proof}
In \cite{Buc}, L. Buchanan completed the classification of prime knots with at most 12 crossings into positive and non-positive knots.
Combining the above inequality $u^{\Delta}(K) \geq c(K)/4$ with the table of values of
$u^{\Delta}(K)$ for prime knots with at most 10 crossings and the values of $a_2(K)$ for prime positive knots with at most 12 crossings,
we obtain Theorem \ref{thm:thmpkdtwo}.
\end{proof}

\begin{restatable}{theorem}{deltatbmin} \label{thm:deltatbmin}
Let $K=C(c_1, c_2, ... , c_m)$ be a two-bridge knot, where $c_1$ is a positive odd integer, and $c_i$ is a positive even integer for $2 \leq i \leq m$, and $m$ is even. 
Among all such knots with a fixed odd crossing number $n \geq 3$, the ${\Delta}$-unknotting number is minimized by $C(n-2, 2)$, and the minimum value is $\frac{1}{2}(n-1)$.
\end{restatable}

\begin{proof} 
By \cite{Naka3}, we have
\[
u^{\Delta}(K)=
\frac{1}{4} \sum_{j=1}^{m/2} (\sum_{i=1}^{j} c_{2i-1}) c_{2j} + \frac{1}{8} ( \sum_{i=1}^{m/2} c_{2i} )^{2} \quad \big(= a_2(K)\big).
\]

Set $m=2k$, and let
\[
S:=\sum_{i=1}^{k}c_{2i}.
\]

Then $S$ is a positive even integer.
Since $K$ is an alternating knot, we have
\[
\sum_{i=1}^{k}c_{2i-1}= n-S.
\]

Since $c_{2j}\geq 2$ for each $j$, we have
\begin{align*}
u^{\Delta}(K)
&=\frac{1}{4}\sum_{j=1}^{k}
\left(\sum_{i=1}^{j}c_{2i-1}\right)c_{2j}
+\frac{1}{8}S^2\\
&\geq \frac{1}{2}\sum_{j=1}^{k}
\left(\sum_{i=1}^{j}c_{2i-1}\right)
+\frac{1}{8}S^2\\
&\geq \frac{1}{2}\sum_{i=1}^{k}c_{2i-1}
+\frac{1}{8}S^2\\
&=\frac{n-S}{2}+\frac{1}{8}S^2\\
&=\frac{n}{2}-\frac{S}{2}+\frac{S^2}{8}\\
&=\frac{n-1}{2}+\frac{(S-2)^2}{8}\\
&\geq \frac{n-1}{2}.
\end{align*}
Therefore,
\[
u^{\Delta}(K)\geq \frac{n-1}{2}.
\]

To show that this lower bound is attained, take
$m=2$, $c_1=n-2$, and $c_2=2$.
Since $n \geq 3$ is odd, $c_1$ is a positive odd integer.
Thus, these values satisfy the assumptions.
In this case,
\begin{align*}
u^{\Delta}(K)
&=\frac{1}{4}(n-2)\cdot 2+\frac{1}{8}\cdot 2^2\\
&=\frac{n-1}{2}.
\end{align*}

\end{proof}

In light of the above results, we propose the following conjecture concerning the $\Delta$-unknotting number for prime positive knots.

\conjmain*



\vspace{1em}






\begin{remark}\label{remark:remarkmonte}
Let $M_n = M(0;(2,1),(3,2),(2n-13,2))$ be a Montesinos knot, where $n$ is a positive even integer with
$n \geq 8$.
Then, we have $u^{\Delta}(M_n) = \frac{1}{2}n$.
\end{remark}

\begin{figure}[htbp]
    \centering
    \includegraphics[width=0.5\linewidth]{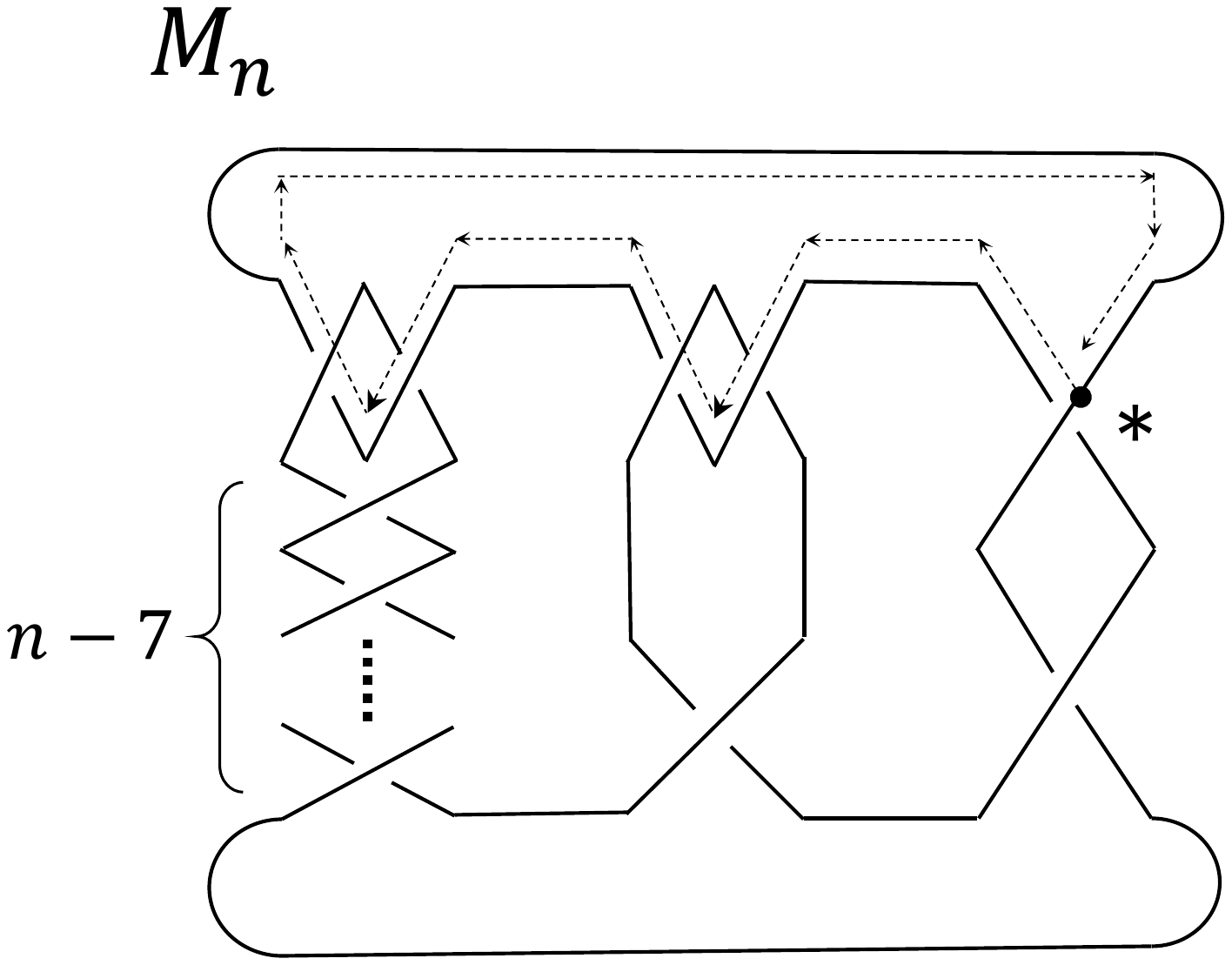}
    \caption{$M_n = M(0;(2,1),(3,2),(2n-13,2))$}
    \label{fig:Mn}
\end{figure}

\begin{proof}
By proceeding in the same way in \cite{Naka3}, we obtain Remark \ref{remark:remarkmonte}.

Let $k_- = M_n$ and $k_+ = C(n-7,2) \# C(1,2)$.

By \cite{Kau},
\[
a_2(C(n-7,2)) \# C(1,2)) - a_2(M_n) = -2.
\]

Thus
\[
a_2(M_n)
= a_2(C(n-7,2) \# C(1,2)) + 2 
= (\frac{1}{2}n - 3) + 1 + 2
= \frac{1}{2}n.
\]

On the other hand, we can deform $M_n$ into $C(n-7,2) \# C(1,2)$ by applying two $\Delta$-moves (see Figure \ref{fig:Mn}).
Therefore, we have
\[
d^{\Delta}(M_n, C(n-7,2) \# C(1,2)) \leq 2.
\]

Furthermore,
\[
u^{\Delta}(C(n-7,2) \# C(1,2)) 
\leq u^{\Delta}(C(n-7,2)) + u^{\Delta}(C(1,2))
= (\frac{1}{2}n - 3) + 1
= \frac{1}{2}n - 2
\]

Thus
\[
u^{\Delta}(M_n)
\leq d^{\Delta}(M_n, C(n-7,2) \# C(1,2)) 
+ u^{\Delta}(C(n-7,2) \# C(1,2)) 
\leq 2 +  (\frac{1}{2}n - 2)
= \frac{1}{2}n
\]

We also have $u^{\Delta}(M_n) \geq |a_2(M_n)|$.
Thus, we conclude that
\[
u^{\Delta}(M_n) = a_2(M_n) = \frac{1}{2}n.
\]
\end{proof}

\begin{remark}\label{remark:remark57}
For prime positive knots with at most 10 crossings,
Conjecture \ref{conj:conjmain} holds.  
\end{remark}

\begin{remark}\label{remark:remark58}
The following table summarizes, for each crossing number, the minimum value of the $\Delta$-unknotting number among prime positive knots and the knots realizing it.
\[
\begin{array}{c|c|c|c|c|c}
c(K) & K & u^{\Delta}(K) & \text{two-bridge knot} & \text{pretzel knot} & \text{Montesinos knot} \\
\hline
3  & 3_1     & 1   & C(1,2) & P(1,1,1) & -- \\
5  & 5_2     & 2   & C(3,2) & P(1,1,3) & -- \\
7  & 7_2     & 3   & C(5,2) & P(1,1,5) & -- \\
8  & 8_{15}  & 4   & --     & --       & M(0; (2,1),(3,2),(3,2)) \\
9  & 9_2, (9_{10}) & 4   & C(7,2) & P(1,1,7) & -- \\
10 & 10_{55}, (10_{154}) & 5   & --     & --       & M(0; (2,1),(3,2),(7,2))\\
\end{array}
\]
There are no prime positive knots with crossing number $4$ or $6$.
\end{remark}

\begin{remark}\label{remark:remark59}
For a prime positive knot $K$, assuming Conjecture \ref{conj:conjmain} holds, we obtain the following lower bounds.
These bounds suggest that the minimum possible value of $u^{\Delta}(K)$ increases as $n$ increases.

\begin{enumerate}
 \item If $n = 3$, then $u^{\Delta}(K) \geq 1$.
 \item If $n = 4$ or $5$, then $u^{\Delta}(K) \geq 2$.
 \item If $n = 6$ or $7$, then $u^{\Delta}(K) \geq 3$.
 \item If $n = 8$ or $9$, then $u^{\Delta}(K) \geq 4$. 
 \item If $n = 10$ or $11$, then $u^{\Delta}(K) \geq 5$.
\end{enumerate}

This leads to the following conjectural classification of prime positive knots:
\begin{itemize}
    \item $u^{\Delta}(K) = 1$ if and only if $K = 3_1 (= C(1,2))$. (see Proposition \ref{thm:thmpkdone})
    \item $u^{\Delta}(K) = 2$ if and only if $K = 5_2 (= C(3,2))$. (see Theorem \ref{thm:thmpkdtwo}(1))
    \item $u^{\Delta}(K) = 3$ if and only if $K = 5_1$ or $7_2 (= C(5,2))$. (see Theorem \ref{thm:thmpkdtwo}(2))
    \item $u^{\Delta}(K) = 4$ if and only if $K = 7_4$, $7_5$, $8_{15} (= M(0; (2,1),(3,2),(3,2)))$, $9_2 (= C(7,2))$, or $9_{10}$.
\end{itemize}

\end{remark}

\begin{remark}\label{remark:remark510}
The following table summarizes the invariants of $C(n-2,2)$ and $M(0;(2,1),(3,2),(2n-13,2))$ (see Theorems \ref{thm:mainv}, \ref{thm:mainw}).
\[
\begin{array}{c|c|c|c}
c(K) & K & u^{\Delta}(K) & u^X(K), g(K), g^*(K), \frac{|\sigma (K)|}{2}\\
\hline
n \text{ odd } (n \geq 3)
& C(n-2,2)=P(1,1,n-2) 
& \frac{1}{2}(n-1) & 1 \\
n \text{ even } (n\geq 8) 
& M(0;(2,1),(3,2),(2n-13,2)) 
& \frac{1}{2} n & 2 \\
\end{array}
\]
In particular, the invariants $u^X(K)$, $g(K)$,
$g^*(K)$, and $\frac{|\sigma(K)|}{2}$ are independent of the crossing number $c(K)$.
\end{remark}

\begin{remark}\label{remark:remark511}
In \cite{Sto}, the following results were established.
\vspace{1em}

Let $K$ be a positive knot with a positive irreducible arborescent diagram of $c$ crossings.
Then $a_2(K) \geq  \frac{1}{2}(c-2)$.

Moreover, the positive twist knot diagrams minimize $a_2(K)$ over all connected irreducible arborescent positive diagrams of odd crossing number.

\vspace{1em}

These results are consistent with Conjecture \ref{conj:conjmain}.
\end{remark}

\vspace{1em}
\section*{Acknowledgments}
The author would like to express his sincere gratitude to Professor Makoto Sakuma for his valuable advice and continuous support.

He also thanks Professors Yasutaka Nakanishi, Yoshiaki Uchida, and Kouki Taniyama for their helpful comments and suggestions.

Finally, he would like to thank his family for their unwavering support and understanding throughout this work.




\end{document}